\documentclass[11pt,psamsfonts,a4paper]{amsart}
\usepackage{amssymb,amsmath,amsthm,amscd}
\chardef\bslash=`\\ 

\makeatletter
\def\verbatim{\interlinepenalty\@M \@verbatim
  \leftskip\@totalleftmargin\advance\leftskip2pc
  \frenchspacing\@vobeyspaces \@xverbatim}
\makeatother
\swapnumbers 

\theoremstyle{plain}
\newtheorem{theorem}{Theorem}[section]

\newtheorem{lemma}[theorem]{Lemma}
\newtheorem{proposition}[theorem]{Proposition}

\theoremstyle{remark}
\newtheorem{remark}[theorem]{Remark}
\newtheorem{example}[theorem]{Example}

\theoremstyle{definition}
\newtheorem{definition}[theorem]{Definition}

\newcommand{\proref}[1]{Proposition~\ref{#1}}

\newcommand{\lemref}[1]{Lemma~\ref{#1}}

\newcommand{\spa}{\operatorname{span}}
\newcommand{\cspa}{\overline{\spa}}

\newcommand{\inv}{^{-1}}

\newcommand{\jay}{ j}
\newcommand{\Gen}{\mathcal G}
\newcommand{\oa}{{\mathcal O}_A}
\newcommand{\uoa}{\widetilde{\mathcal O}_A}
\newcommand{\ome}{\Omega_A}
\newcommand{\ometil}{\widetilde{\Omega}_A}
\newcommand{\coker}{\operatorname{coker}}
\newcommand{\ima}{\operatorname{Im}}
\newcommand{\odd}{\operatorname{odd}}
\newcommand{\even}{\operatorname{even}}

\def\urowa{\tilde{R}_A}
\def\rowa{R_A}
\def\zrowa{{\mathfrak R}_A}
\def\zurowa{\tilde{\mathfrak R}_A}
\def\idx{{\operatorname{id}}_x}

\def\zerocol{\mathbf 0}
\def\onecol{\mathbf 1}
\def\bfone{\mathbf I}
\def\N{\mathbb N}

\begin{document}

  \title[$K$-theory of $\oa$ for infinite $A$]
  {The $K$-theory of Cuntz-Krieger algebras \\for infinite matrices}
\date{March 8, 1999}

\author[R. Exel]{Ruy Exel${}^*$}
\thanks{${}^*$) Partially supported by CNPq.}
\address{Departamento de Matem\'{a}tica, Universidade Federal de Santa
Catarina,\  88010---970 Florian\'{o}polis, 
Santa Catarina, Brazil.}
\email{exel@mtm.ufsc.br}
\author[M. Laca]{Marcelo Laca${}^\dagger$}
\thanks{${}^\dagger$) Supported by the Australian Research Council.}
\address{Department of Mathematics, The
        University of Newcastle,  NSW  2308, Australia.}
\email{marcelo@math.newcastle.edu.au}

\begin{abstract} 
We compute the $K$-theory groups of the Cuntz-Krieger $C^*$-algebra $\oa$
associated 
 to an infinite matrix $A$ of zeros and ones.
\end{abstract}

\maketitle

\section*{Introduction} 
    Let $\oa$ be the Cuntz-Krieger $C^*$-algebra associated in \cite{cun-kri} 
    to a $0$-$1$ matrix $A$.
   Cuntz  showed in \cite{cun2} that  $K_1(\oa)$ and  $K_0(\oa)$ are,
   respectively, the kernel and co-kernel of the transformation 
              $$
              (I-A^t) : {\mathbb Z}^n \to {\mathbb Z}^n.
              $$
              
  In the case of a countably infinite matrix $A$ with only finitely many `ones'
  in each row, this result has been extended to the
  corresponding Cuntz-Krieger algebras $\oa$ in \cite{para}, with
  ${\mathbb Z}^n$ replaced by the countably infinite direct sum of 
  copies of $\mathbb Z$, see also \cite{kpr}.

  The definition of $\oa$ has been extended to general infinite 
$0$-$1$ matrices $A$ in \cite{infinoa}. However, when $A$ is not row-finite,
$(I-A^t)$ does not map the infinite direct sum of copies of $\mathbb Z$ into itself
  anymore, so the obvious analogue of the above characterization of $K_0(\oa)$
  does not make sense. 
 The purpose of this paper is to show that also for infinite $A$ it is the case that $K_1(\oa)$ and
  $K_0(\oa)$ are the kernel and co-kernel of 
  $(I-A^t)$, once the appropriate domain and co-domain are established.
  It turns out that the right domain is simply the direct sum over $\Gen$ of copies of 
  $\mathbb Z$, but the
  construction of the right co-domain is a rather subtle point that involves
  the matrix $A$ in a nontrivial way. Indeed, as we will see in 
  Section \ref{themap}, this co-domain turns out to be the ring
  of functions on $\Gen$ generated by the delta functions and the rows of $A$, 
  viewed as functions on $\Gen$.
  
  As a byproduct of studying this ring and its associated $C^*$-algebra we also obtain
   a new formulation of the relations defining $\oa$, \proref{rels-revisit},
   which exhibits the relevant underlying algebraic structure in an explicit manner.
 
\section{Preliminaries.}\label{prelim}

  Let $\Gen$ be a set and suppose 
  $A =  \{ A(i,j) \}_{i,j\in\Gen}$ is a $0$-$1$ matrix
  which we assume has no identically zero rows.  
  According to \cite[\S 7]{infinoa}, the {\em unital  Cuntz-Krieger $C^*$-algebra $\uoa$} 
   is the universal unital $C^*$-algebra generated by 
  partial isometries $\{S_i: i\in \Gen\}$ such that
   for $i,j \in \Gen$,
 \begin{enumerate}
  \item[(i)]  $S_i^*S_i$ and $S_j^*S_j$ commute,
  \item[(ii)]  $(S_iS_i^*)(S_jS_j^*) = \delta_{i,j} S_i S_i^*$,
  \item[(iii)]  $(S_i^*S_i)( S_jS_j^* )= A(i,j) S_jS_j^*$, and
  
  \smallskip
   \item[(iv)]  $ 
   \prod_{x\in X} S_x^* S_x \prod_{y\in Y} (\bfone - S_y^* S_y)\  =\ 
     \sum_{j\in\Gen} A(X,Y,j) S_j S_j^*
     $
     whenever $X$ and $Y$ are finite subsets of $\Gen$ such that the function 
     $
    j\in \Gen \mapsto  A(X,Y,j):= \prod_{x\in X} A(x,j) \prod_{y\in Y} (1 - A(y,j))
     $
     is finitely supported. 
    \smallskip
 \end{enumerate}
  The {\em Cuntz-Krieger algebra $\oa$} is defined to be the $C^*$-subalgebra
  of $\uoa$ generated by the $S_i$.
  In fact,  $\oa$ itself can be viewed as the universal (not necessarily unital)
  $C^*$-algebra with the above presentation, provided that we interpret
  the appearance of the unit $\bfone$ in (iv) as follows:
     if $\bfone$ still appears after expanding the product in the left hand side of (iv),
     i.e. if there is a finite subset $Y$ such that $A(\emptyset, Y,j)$
     is finitely supported, then the resulting condition implicitly states that 
     $\oa$ is unital, in which case it is equal to $\uoa$.
     When $\oa$ is not unital, it is an ideal of codimension $1$ in $\uoa$, so 
   that $\uoa$ is indeed the unitization of $\oa$.
  
  Let $\mathbb F$ be the free group on $\Gen$. A key step in
\cite{infinoa} is
  to realize  $\uoa $ and $\oa$
  as  crossed products by partial actions of $\mathbb F$ on commutative algebras.
  Since  the relevant  partial dynamical systems
  $(C(\ometil), {\mathbb F}, \alpha)$
  and 
  $(C_0(\ome), {\mathbb F}, \alpha)$
  are central to our discussion, we summarize the basic definitions 
   and results below. The reader is referred to \cite{infinoa} 
   for the details.
  
  Denote by $2^{\mathbb F}$ the topological space $\{0,1\}^{\Gen}$ equipped with the
  product topology. Whenever convenient we will identify $2^{\mathbb F}$ with the space of
  all subsets of $\mathbb F$ in the usual way.  As in 
  \cite[\S 5]{infinoa}, denote by $\ome^\tau $ the subset of $2^{\mathbb F}$ given by 
  $$
  \left\{ 
  \xi \in 2^{\mathbb F}: \begin{array}{rl}  
   &e\in \xi, \ \xi \text{ is convex, } \\
     &\text{ if }\omega \in \xi
   \text{ then there is at most one } y\in {\Gen} \text{ with }\omega y \in \xi, \\
     & \text{ if }\omega, \omega y \in \xi  \text{ then }
  ( \omega x\inv \in \xi \iff A(x,y) = 1).
  \end{array}
  \right\}.
  $$ 
  We recall that a subset $\xi$ of $\mathbb F$
  is convex if it contains the shortest path between any two of its elements \cite[Definition 4.4]{infinoa}. 
  It is helpful to notice that
 the set $\xi$ is convex and contains $e$ if and only if it contains all the
initial subwords of each
   one of its elements.

  Let $\mathbb F^+$ be the 
  unital semigroup generated by $\Gen$ in $\mathbb F$;
  we say that an element $\xi$ is {\em unbounded} if it has no upper bound in $\mathbb F$
  with respect to the left order defined by $s \leq t$ when $s\inv t \in \mathbb F^+$.
  We define $\ometil$ to be the closure of the set of unbounded elements in $ \ome^\tau $
   and we let
  $\ome := \ometil \setminus \{\{ e \}\}$. The singleton $\{e\}$ may well
  not be in $\ometil$, in which case $\ome = \ometil$; otherwise $\ometil $ is
  the one-point compactification of $\ome$.
  There is a partial action of $\mathbb F$ on
  $\ometil$ by partial homeomorphisms $h_t : \Delta_{t\inv} \to \Delta_t$
  where the range $\Delta_t$ of $h_t$ is the clopen set
  $\{ \xi \in \ometil : t\in \xi\}$, and $h_t$ is defined by $ h_t (\xi) = t\xi$ for $\xi \in
\Delta_{t\inv}$.
   At the level of $C^*$-algebras, the partial action is given by
   the partial automorphisms $\alpha_t$ of $ D_{t\inv}:=
C_0(\Delta_{t\inv})$  to
   $D_t := C_0(\Delta_t)$ given by $\alpha_t(f)(\xi) = f(t\inv \xi)$
   for $f\in D_{t\inv} $ and $\xi \in \Delta_t$.
   
   By Theorems 7.10 and 8.4 of \cite{infinoa}, there are canonical isomorphisms
   $$
   \uoa \cong C(\ometil) \rtimes \mathbb F \qquad \text{ and  }\qquad
    \oa \cong C_0(\ome) \rtimes \mathbb F,
   $$
  in which, for each $x\in \Gen$, the generating
   partial isometry $S_x$ is mapped to the partial isometry associated to  
   $x \in \Gen \subset \mathbb F$ in the crossed product. 
  Under the  corresponding identification of $C(\ometil)$ to a subalgebra of $\uoa$
  the projection $S_xS_x^*$ is
  identified to  the characteristic function of the set $\Delta_x$ and the
  projection $S_x^*S_x$ is identified to the
  characteristic function of $\Delta_{x\inv}$.     
  
\section{$K$-theory for partial actions of free groups.} \label{ktheory-sec}
 
  In this section we obtain a generalization,
  to the case of infinitely many generators, of McClanahan's exact sequence for
  crossed products by partial actions of free groups \cite{mac}, cf. \cite{ex-pa}.
  Since we need to distinguish the free groups on various sets of generators, we
  use ${\mathbb F}_{\mathcal S}$ to denote the free group on a set $\mathcal S$.
  
  Suppose $\alpha$ 
  is a partial action of ${\mathbb F}_{\Gen}$ on a $C^*$-algebra
  ${\mathcal A}$. For each $t\in{\mathbb F}_\Gen$  denote by $D_t$ the range of the
  partial automorphism $\alpha_t$. By Theorem 6.2 of \cite{mac} we have the following
    exact sequence, valid for finite $\Gen $:
\stepcounter{theorem}\begin{equation} 
\label{pvseq}
\begin{CD}
\bigoplus_{x\in \Gen} K_0(D_x)   @>{\ \sum_{x\in {\Gen}}(\idx - \alpha_x\inv)_*\ }>>
K_0({\mathcal A}) @>{\ \ \  (i_A)_*\ \ \  }>> K_0({\mathcal A}\rtimes_{\alpha}{\mathbb F}_{\Gen}) \\
 @AAA @. @VVV  \\  
K_1({\mathcal A}\rtimes_{\alpha}{\mathbb F}_{\Gen})  @<{\ \ \  (i_A)_*\ \ \ }<<
  K_1({\mathcal A}) 
  @<{\ \sum_{x\in {\Gen}}(\idx - \alpha_x\inv)_*\ }<<
\bigoplus_{x\in {\Gen}} K_1(D_x),
\end{CD}
\end{equation}

\medskip\noindent
   in which $\idx$ is the embedding of $D_x$ in $\mathcal A$ and 
   $i_A$  is the embedding of
  ${\mathcal A}$ in ${\mathcal A}\rtimes_{\alpha} {\mathbb F}_{\Gen}$.

  We will need this exact sequence for {\em infinite} $\Gen$,
  and we derive it from the finite case by first realizing
  ${\mathcal A}\rtimes_{\alpha}{\mathbb F}_{\Gen}$ as a direct limit and
  then using the continuity of $K$-theory.
  
 \begin{proposition} 
 \label{pvseq-lim}
  Let $\Gen$ be an infinite set.
  Suppose $\alpha$ is a partial action of ${\mathbb F}_{\Gen}$ on the
  $C^*$-algebra $\mathcal A$ and let $\Lambda$ denote the collection of finite subsets of
  $\Gen$, directed by inclusion. For every $\lambda \in \Lambda$
  denote also by $\alpha$ the partial action restricted to the subgroup 
  ${\mathbb F}_\lambda$. Then there is a directed system 
  $$
  \Phi_\lambda^{\mu}:
  {\mathcal A}\rtimes_{\alpha}{\mathbb F}_{\lambda} \to 
  {\mathcal A}\rtimes_{\alpha}{\mathbb F}_{\mu}
  , \qquad \lambda \subset \mu
  $$ 
  determined by the inclusion ${\mathbb F}_{\lambda} \subset {\mathbb F}_{\mu}$,
  and 
  $
  {\mathcal A}\rtimes_{\alpha}{\mathbb F}_\Gen = 
  \lim_\lambda {\mathcal A}\rtimes_{\alpha}{\mathbb F}_\lambda.
  $
 
 \end{proposition}

  \begin{proof}
  That such a directed system exists follows easily from
  the universal property of crossed products with respect to covariant
  representations of the corresponding partial dynamical systems 
  \cite[Theorem 1.4]{elq}, cf. \cite[Propositions 2.7 and 2.8]{mac}.
  It also follows that for each $\lambda$ there is a homomorphism $\Phi_\lambda^\Gen$ 
  of ${\mathcal A}\rtimes_{\alpha}{\mathbb F}_\lambda$ into the
  crossed product ${\mathcal A}\rtimes_{\alpha}{\mathbb F}_\Gen$ such that
  $\Phi_\lambda^\Gen = \Phi_\mu^\Gen \circ \Phi_\lambda^{\mu} $.
  
  In order to prove that ${\mathcal A}\rtimes_{\alpha}{\mathbb F}_\Gen$ 
  is the direct limit of the directed system it suffices to prove that if 
  $$\psi_\lambda : {\mathcal A}\rtimes_\alpha {\mathbb F}_\lambda \to B $$ 
  is a family of homomorphisms into a $C^*$-algebra $B$ that is compatible with the directed system in
  the sense that $ \psi_\lambda = \psi_\mu \circ \Phi_\lambda^{\mu}$, 
  then there exists
  a unique homomorphism $\sigma$ of ${\mathcal A}\rtimes_{\alpha}{\mathbb F}_\Gen$
to $B$
  such that   the following diagram commutes
 $$ 
\begin{CD}
{\mathcal A}\rtimes_{\alpha}{\mathbb F}_\lambda  @>{ \Phi_\lambda^\Gen }>>  
{\mathcal A}\rtimes_{\alpha}{\mathbb F}_\Gen\\
 @. @VV{\sigma}V  \\  
  \ \   @.   B
\end{CD}
\hspace{-8em}
\begin{picture}(35,35)
               \put(10,8){\vector(3,-2){35}}
               \put(25,0){$\psi_\lambda$}
\end{picture}
$$

  Assume the $\psi_\lambda$ satisfy  $\psi_\mu \circ \Phi_\lambda^{\mu} = \psi_\lambda$ and suppose,
  without loss of generality, that
  $B$ is contained in $B(H)$ for some Hilbert space $H$.
  Thus each $\psi_\lambda$ is of the form $\pi_\lambda \times V_\lambda$, for a 
  covariant representation $(\pi_\lambda, V_\lambda)$
  of the system $({\mathcal A},{\mathbb F}_{\lambda},\alpha)$.
  
  Since $\psi_\lambda =\psi_\mu \circ \Phi_\lambda^{\mu} $, the representation 
  $\pi_\lambda$ of $\mathcal A$ does not depend on $\lambda$. Suppose  $x \in \lambda \subset \mu$
   and let $v_x$ be the partial isometry corresponding to $x$ 
  in $\mathcal A \rtimes_\alpha \mathbb F_\lambda$,   then 
  $V_\lambda (x) = \psi_\lambda (v_x)  = (\psi_\mu \circ \Phi_\lambda^\mu ) (v_x)$. But
  $\Phi_\lambda^\mu (v_x)$ is the partial isometry  corresponding to $x$
  in $\mathcal A \rtimes_\alpha \mathbb F_\mu$, so we conclude that $V_\lambda (x) = V_\mu(x)$. 
  Thus the $V_\lambda$ determine 
  a partial representation $V$ of ${\mathbb F}_{\Gen}$ such that $(\pi,V)$ is covariant.
  By \cite[Theorem 1.4]{elq}, there exists a representation $\pi \times V$ of ${\mathcal
  A}\rtimes_{\alpha}{\mathbb F}_\Gen$ in $B$ such that 
  $$\psi_\lambda = (\pi \times V) \circ \Phi_\lambda^\Gen.$$
  This proves that ${\mathcal A}\rtimes_{\alpha}{\mathbb F}_\Gen$ has the property
  that characterizes the inductive limit up to canonical isomorphism, finishing the proof.
  \end{proof} 

  \begin{proposition}
  The sequence \eqref{pvseq} is exact also for $\Gen$ infinite.
  \end{proposition}
  \begin{proof}
  By continuity of $K_*$, every term of the sequence is
  the limit of its finite counterparts. Since the connecting 
  maps are compatible with these limits, the sequence is exact.
  \end{proof} 
  
  \begin{proposition} 
  Define a map
  $L_\alpha : \bigoplus_{x\in \Gen}C(\Delta_x, {\mathbb Z}) \to C(\ometil, {\mathbb Z})$ by
  $$
  L_\alpha f := \sum_{x \in \Gen}(\idx - \alpha_x\inv) f_x 
        \qquad \text{ for } 
        f = ( f_x)_{x\in \Gen} \in \bigoplus_{x\in \Gen}C(\Delta_x, {\mathbb Z})
  $$
  Then
  \begin{enumerate}
  
  \item[(i)]   $K_0 (\uoa) = C(\ometil,{\mathbb Z})/ \ima L_\alpha$,              
  
  \item[(ii)]  $K_0 (\oa) = C_0(\ome,{\mathbb Z})/ \ima L_\alpha $, and
  
  \item[(iii)] $K_1(\uoa) =  K_1(\oa) = \ker L_\alpha $.

  \end{enumerate}
  
  \end{proposition}
  \begin{proof}
  The sets  $\ometil$ and $\Delta_x$ for $x\in \Gen$ are totally disconnected, so
  $$ 
  K_1(C(\ometil)) = 0 =  K_1(C(\Delta_x)),
  $$
while
  $$
  K_0(C(\ometil)) = C(\ometil, {\mathbb Z}) \qquad \text{ and }\qquad
  K_0 (C(\Delta_x)) = C(\Delta_x, {\mathbb Z}) .
  $$
  With these simplifications
  the exact sequence \eqref{pvseq}, when applied to the crossed product 
  $\uoa = C(\ometil) \rtimes_\alpha {\mathbb F}_{\Gen}$, gives
 $$ 
 \begin{CD}
   \bigoplus_{x\in \Gen} C(\Delta_x, {\mathbb Z})  @>{\ 
   L_\alpha\ }>> C(\ometil, {\mathbb Z}) @>{\ \ \ \ \ }>>  K_0(\uoa) \\
   @AAA @. @VVV  \\  
    K_1(\uoa)  @<{\ \ \ \ \ }<<
    0  @<{\ \ }<< 0,
 \end{CD}
 $$ 
 
 \medskip\noindent
  {}from which the expressions for $K_0(\uoa)$ and $K_1(\uoa)$ follow easily.
  The expression for $K_0(\oa)$ is obtained via a similar argument on the
  crossed product $\oa \cong C(\ome) \rtimes {\mathbb F}_{\Gen}$ by the partial
  action restricted to $\ome$.
  \end{proof} 
  
\section{The $C^*$-subalgebra generated by the projections $P_i$ and $Q_i$.}  
   
   The expressions of $K_*(\uoa)$ and $K_*(\oa)$ obtained in the preceding section
   are of limited usefulness because of the difficulties associated 
   with computing the kernel and co-kernel of the map $L_\alpha$ in an explicit form. 
   In order to obtain computable expressions that involve the matrix $A$ in a 
   more direct way we will need to understand the structure of the $C^*$-subalgebra
   of $\oa$ generated by the projections $P_i$ and $Q_i$.

   As before, $A$ will denote a $0$-$1$ matrix over a set $\Gen$.
   For each $i\in\Gen$  denote by $\rho_i$ the $i^{th}$ row of the matrix $A$, 
   viewed as the function from $\Gen$ to $\{0,1\}$, defined by $\rho_i(j) = A(i,j)$
   for $j\in \Gen$. As usual, denote  by $\delta_i$ 
   the $i^{th}$ row of the  identity  matrix over $\Gen$.
   \begin{definition}
   \label{RowaDef}
   The $C^*$-subalgebra of  $\ell^{\infty} ({\Gen})$ generated by the rows of $A$ and of $I$  
   will be denoted by $\rowa$:
    $$
     \rowa := C^*(\{\rho_i ,\delta_i: i\in
    {\Gen}\}) 
    \subset \ell^{\infty} (\Gen).
    $$
    We will also work with the unitization of $\rowa$, which is
    defined by
    $$
    \urowa := C^*(\rowa \cup \{\onecol\}) \subset \ell^{\infty} (\Gen).
    $$
    \end{definition}
  
  The first step in establishing the relation between $\rowa$ and $\oa$
  is to obtain a concrete picture of the spectrum of $\rowa$.

  Denote by $\tilde{\Gen}$ the one-point compactification
  of the (discrete) space $\Gen$,
  so that $\tilde{\Gen} = \Gen \cup \{\star\}$ 
  if $\Gen $ is infinite, and $\tilde{\Gen} = \Gen $ if it is finite.
  
  For each $j\in {\Gen}$ let $c_j$ be the $j^{th}$ column of the matrix $A$, viewed as
  the function on $\Gen$ defined by
  $c_j(i) = A(i,j)$ for $i \in \Gen$. Thus $c_j $ is 
  a $\{0,1\}$-valued function on $\Gen$, i.e. an element of $\{0,1\}^\Gen$, and, as such, it may
  also be viewed as a subset of $\Gen$. 
  There is no harm in shifting from one point of view to the other
  because  $\{0,1\}^\Gen$ is homeomorphic to the power set of $\Gen$.
 
 \begin{definition}
   Let $\tilde \Gamma_A$ be the compact space obtained as 
  the closure in $ \tilde{\Gen} \times \{0,1\}^\Gen$ 
  of the set of elements of the form $(j, c_j)$ for $j \in \Gen$: 
  $$
  \tilde{\Gamma}_A : = 
  \overline{\{ (j, c_j) \in \tilde{\Gen} \times \{0,1\}^\Gen: j \in {\Gen}\}}.
  $$
  The map  $j \mapsto (j, c_j)$ embeds $\Gen$ in $\tilde{\Gamma}_A$ 
  as a relatively discrete dense subset,
    (the graph of the map $j \mapsto c_j$), so $\tilde{\Gamma}_A$ may
   be seen as a compactification  
  of $\Gen$ that reflects the column structure of the matrix $A$.
  
  We will also consider the set
  $$
  \Gamma_ A : = \tilde{\Gamma}_A  \setminus \{(\star,\zerocol)\},
  $$ 
  with $\zerocol$ denoting the identically zero function on $\Gen$.
  \end{definition}
  
  Notice that
   the point $(\star,\zerocol)$ may actually not belong to $\tilde\Gamma_A$, in which case
   $\tilde{\Gamma}_A = \Gamma_A$.  When $(\star, \zerocol)$ is in 
   $\tilde \Gamma_A$, our notation agrees with the standard one,
   in that $\tilde\Gamma_A$ is the one-point
   compactification of the locally compact space $\Gamma_A$. 
   
  The points of $\Gamma_A$ are either of the form $(j, c_j)$, 
  with $j\in \Gen$, or
  $(\star, c)$, with $\star$ the point at infinity of $\tilde{\Gen}$ and 
  $c $ an accumulation point, in 
  $ \{0,1\}^\Gen \setminus \{\zerocol\}$, of the collection of columns 
  $\{c_j\}_{j\in {\Gen}}$. We stress that the collection of columns appears here as an indexed set:
  specifically,  $c\in \{0,1\}^\Gen$ is an accumulation point of the collection of columns
  if, for every
  neighborhood $V$ of $c$ in $\{0,1\}^\Gen$, the set 
  $\{ j \in \Gen : c_j \in V\}$ is infinite;  equivalently,  if
  for every finite  $F \subset \Gen$ there are infinitely many $j \in \Gen$ 
  for which $c_j$ is equal to $c$ on the subset $F$. 

  \begin{remark}
  In the finite and row-finite situations the rows already are in the $C^*$-algebra generated
  by the $\delta_i$'s, so $\rowa$ and $\Gamma_A$ are easy to characterize:
  \begin{itemize}
  \item If  $\Gen$ is finite, then it is already compact, so $\Gamma_A \cong
  \Gen$, and $\rowa \cong {\mathbb C}^n$, where $n$ is the number of elements in $\Gen$.

  \item If $\Gen$ is infinite but the matrix $A$ is row-finite, then the identically
  zero column is the unique accumulation point of $\{c_j\}_{j\in \Gen}$, so
  $\tilde{\Gamma}_A = \{ (j,c_j): j\in \Gen\} \cup \{(\star, \zerocol)\}$ and
  $\Gamma_A \cong \Gen$. In this case $\rowa $ is $C^*(\{\delta_i: i\in \Gen\}) = C_0(\Gen)$.
\end{itemize}
  \end{remark}

  \begin{proposition}
  \label{SpectrumRowa}
  View $\Gen $ as a dense subset of $\tilde \Gamma_A$ by 
  identifying $j$ to $(j,c_j)$. 
  Then every function $f\in \urowa$ has a unique extension $\hat f \in C(\tilde\Gamma_A)$,
  and the map $f\mapsto \hat f$ gives isomorphisms
  $$
  \urowa \cong C(\tilde\Gamma_A) 
  \qquad \text{ and }\qquad
  \rowa \cong C_0(\Gamma_A). 
  $$ 
   \end{proposition}
   \begin{proof} 
   Denote a typical element of $\tilde\Gamma_A $ by $(x,c)$ so that $x \in \Gen \cup \{\star\}$
   and $c $ is either a column or an accumulation point of columns. We identify
   $x\in \Gen$ to the point $(x,c_x)$, and so we may view the generators  $\delta_i$ and $\rho_i$
   of $\urowa$, as functions on a dense subset of
   $\tilde\Gamma_A$. Clearly it suffices to extend these generators to 
   functions $\hat\rho_i$ and $\hat\delta_i$ on $\ometil$;  
   uniqueness of the extension is a consequence of the density of $\Gen$ in $\tilde \Gamma_A$.
   
   For each $i\in \Gen$ the function defined by $\hat\delta_i(x,c) = \delta_{i,x}$
   is continuous and extends $\delta_i$ to $\tilde \Gamma_A$,  
   and the function defined by $\hat \rho_i (x,c) = c(i)$ is  continuous and extends
   $\rho_i$ to  $\tilde \Gamma_A$. Hence we may view $\urowa$ as a 
   unital subalgebra of $C(\tilde\Gamma_A)$; to prove equality 
   it suffices to show that $\urowa$ separates points of $\tilde\Gamma_A$. 
  It is clear that for each $i\in \Gen$ the function $\hat\delta_i$ separates 
  the point $(i,c_i)$ from all
  other points of $\tilde\Gamma_A$. It remains to consider the case 
  in which the two points to be separated are
  $(\star,c) $ and $(\star,c')$, with $c$ and $c'$ different accumulation points of the
  collection of columns.
  Since $c \neq c'$ there exists $i\in \Gen$ such that 
  $\hat\rho_i (\star,c) = c(i) \neq c'(i) = \hat\rho_i(\star,c')$, so the points are separated by the 
  unique extension of $\rho_i$. This yields the first isomorphism.
  
  Assume now that the point $(\star,\zerocol)$ is in
  $\tilde\Gamma_A$. Then
  $\rowa$ is an ideal of codimension $1$ in $\urowa$, all of whose generators
  vanish at $(\star, \zerocol)$. Since
   $C_0(\Gamma_A)$ is the ideal of codimension $1$ of $C(\tilde\Gamma_A)$
  consisting of all functions vanishing at $(\star,\zerocol)$,
  the isomorphism of $\urowa$ to $C(\tilde\Gamma_A)$ carries $\rowa$ onto $C_0(\Gamma_A)$.
   \end{proof}
 
  The next step is to relate $\Gamma_A$ to the generalized Cuntz-Krieger 
  relations (i) --- (iv) of Section \ref{prelim}.
  Recall, from Definitions 5.5 and 5.6 of \cite{infinoa}, that every point $\xi$ in the
  spectrum $\ometil$ of the unital Cuntz-Krieger relations has 
  a {\em root at the identity}
  $R_\xi(e) := (\xi\inv \cap \Gen)$, and a 
  {\em stem} $\sigma(\xi) := \xi \cap {\mathbb F}^+$, which is a
  positive word. Here a (finite or infinite) positive word 
  is represented by the collection of its
  finite initial subwords, cf. Definition 5.2 and Proposition 5.4 of \cite{infinoa}.

  \begin{proposition} 
  \label{psihat}
  For each $\xi \in \ometil$ let  $R_\xi(e)$ be the root of $\xi$ at the identity, viewed as an
  element of $2^\Gen$,  and let $\sigma(\xi)_1$ denote the initial letter
  of the stem  $\sigma(\xi)$ of $\xi$, with the convention that if the stem is trivial,
  then $\sigma(\xi)_1 = \star \in \tilde \Gen$.
  
  \begin{enumerate}
  \item[(i)] The map defined by 
 $$
   \xi \in \ometil \mapsto 
  (\sigma(\xi)_1 , R_\xi(e))
 $$ 
 is a continuous surjective map from $\ometil$ onto $\tilde{\Gamma}_A$,
 mapping $\ome$ onto ${\Gamma}_A$.
 
 \item[(ii)] The homomorphism
  $
  \psi: \urowa \to C(\ometil)
  $ defined by $\psi(f)(\xi) = \hat f (\sigma(\xi)_1 , R_\xi(e))$
  is  injective  and satisfies 
  $$
  \psi(\delta_i) = P_i \qquad \text {and } \qquad \psi(\rho_i )= Q_i,
  $$
giving  isomorphisms 
  $$
  \urowa \cong C^*(\{\bfone, P_i,Q_i: i\in \Gen\}) 
  \qquad \text{ and } \qquad
  \rowa \cong C^*(\{ P_i,Q_i: i\in \Gen\}).
  $$
  \end{enumerate}
  \end{proposition}
  
  \begin{proof}
  First we prove continuity of the map defined in (i).
   The map $\xi \in 2^{\mathbb F} \mapsto \xi \cap \Gen \in 2^\Gen$
   is obviously continuous in the respective product topologies.
   The set $\xi \cap \Gen $
   is either empty or equal to the singleton $\{\sigma(\xi)_1\}$ 
   with $\sigma(\xi)_1$ the initial letter of the stem of $\xi$.
   Since the  map $x \mapsto \{x\}$ and 
   $\star \mapsto \emptyset$ 
   establishes a homeomorphism of $\tilde{\Gen}$ to the subspace 
   $\{c\in 2^\Gen : \#(c) \leq 1\}$ of $2^\Gen$,
   it follows that the first coordinate of $\psi$ is continuous.   
   A similar argument on the map
   $\xi \in 2^{\mathbb F} \mapsto R_\xi(e) \in 2^\Gen$ shows that the second
   coordinate is continuous too.
    
   To prove surjectivity it suffices to show that the image of $\ometil$ contains every
   $(x,c_x)$ with $x\in \Gen$, because these elements form a dense set and the image of the map
   is closed.  Let $(x,c_x)$ be such an element and take $\xi$ in $\ometil$ such that $\xi \cap \Gen =
   \{x\}$  
   (e.g. take $\xi$ to be the unbounded element associated as in \cite[Proposition 5.13]{infinoa}
    to an infinite word starting with $x$, which exists because
   $A$ has no identically zero rows).
   Since $x\in \xi\cap \Gen$ and $\xi$ is in $\ometil$, we have that $ R_\xi(e)= c_x$, and
   hence 
   $(\sigma(\xi)_1 , R_\xi(e)) = (x, c_x) $. 
   Notice that if $ (\star,\zerocol) $ is in $ \tilde\Gamma_A$, then $\{e\}$ is in
   $\ometil$ by \cite[Proposition 8.5]{infinoa},
   and clearly $\{e\}$ is the only point of $\ometil$ mapped onto $(\star,\zerocol)$.
   Hence $\ome = \ometil \setminus \{\{e\}\}$ is mapped onto 
   $\Gamma_A = \tilde\Gamma_A \setminus \{(\star,\zerocol)\}$.
 
  Next we prove part (ii).  
  The  homomorphism $\psi$ is injective, because it is the adjoint map of
  the surjective continuous map constructed in (i) 
  composed with the identification of $\urowa$ to  $C(\tilde\Gamma_A)$
  given in \proref{SpectrumRowa}.
  To prove that $\psi(\delta_i) = P_i$ and $\psi(\rho_i) = Q_i$ it suffices to
  check equality on the dense subset of unbounded elements of $\ometil$.
  Let $\xi$ be an unbounded element.
  Then
  $$
  \psi(\delta_i) (\xi) = \hat{\delta}_i(\sigma(\xi)_1 , R_\xi(e)) =
  \delta_i(\sigma(\xi)_1)
  $$
  so $\psi(\delta_i)$ is the characteristic function of the set
  $\Delta_{i} : = \{\xi \in \ometil: \sigma(\xi)_1 = i\}$,
  which is precisely $P_i$.
  As for $\rho_i$, we have that 
  $$
  \psi(\rho_i) (\xi)= \hat{\rho}_i(\sigma(\xi)_1 , R_\xi(e)) =  R_\xi(e)(i) 
  $$
  so  $\psi(\rho_i)$ is the characteristic function of 
  the set  $\Delta_{i\inv} : = \{\xi \in \ometil: i\inv \in \xi\}$, i.e.,
   $Q_i$.   
  \end{proof} 

Before launching into the $K$-theory computations, we digress momentarily to 
take a new look at the relations defining of $\oa$ for infinite $A$. The 
material of the remainder of this section 
is not needed for the main results presented in the next section;
it is included here because it may help clarify the role of 
relation (iv) in the definition of $\oa$.
We begin with a characterization of $\rowa$ in terms of generators and relations.

\begin{proposition}
\label{rowa-univ}
 Suppose $\{P_i,Q_i: i \in \Gen\}$ is  a family of projections satisfying the relations
 \begin{enumerate}
  \item[(i)]  $Q_i$ and $Q_j$ commute,
  \item[(ii)]  $P_i P_j = \delta_{i,j} P_j$,
  \item[(iii)]  $Q_i P_j= A(i,j) P_j$, and
  
  \smallskip
   \item[(iv)]  $ 
   \prod_{x\in X} Q_x \prod_{y\in Y} (\bfone - Q_y)\  =\ 
     \sum_{j\in\Gen} A(X,Y,j) P_j
     $
     whenever $X$ and $Y$ are finite subsets of $\Gen$ such that the function 
     $
    j\in \Gen \mapsto  A(X,Y,j):= \prod_{x\in X} A(x,j) \prod_{y\in Y} (1 - A(y,j))
     $
     is finitely supported.
    \smallskip
 \end{enumerate}
 Then the maps $\delta_i \mapsto P_i$, and $\rho_i \mapsto Q_i$ extend to a $C^*$-algebra homomorphism
 of $\rowa$ onto $C^*(\{P_i,Q_i: i\in \Gen\})$  and of $\urowa$ onto 
 $C^*(\{\bfone, P_i,Q_i: i\in \Gen\})$.
\end{proposition}

\begin{proof}
We prove the unital case first.
Let $C^*(\{p_i,q_i: i\in \Gen\})$ be the universal unital $C^*$-algebra of the above relations.
The $p_i$ are mutually orthogonal by (i)  and,
on taking adjoints, (iii) 
implies that the $p_i$ commute with the $q_j$. Hence,
$C^*(\{p_i,q_i: i\in \Gen\})$ is abelian and its spectrum, which we denote by $X$,
is totally disconnected. Since the generators of $\urowa$
clearly satisfy the relations,  the universal property gives a surjective homomorphism
of $C^*(\{p_i,q_i: i\in \Gen\}) = C(X)$ onto
$\urowa= C(\tilde\Gamma_A)$, which is adjoint to an injective continuous map of $\tilde\Gamma_A $ into
$X$. We need to show that this map is surjective.

Let $x\in X$. 
Since the $p_i$ are pairwise orthogonal projections there exists 
at most one $i\in \Gen$ such that $p_i(x) \neq 0$. 

Suppose first that there exists such an element $i_0$, so that
$p_i(x) = \delta_{i,i_0}$. We claim that
$x$ is the image of $(i_0, c_{i_0}) \in \tilde \Gamma_A$. To prove this it suffices to
check that $p_i(x) = \hat\delta_i(i_0, c_{i_0})$ and that 
$q_i(x) = \hat \rho_i (i_0, c_{i_0})$. The first equality is easily verified,
since both sides are equal to $\delta_{i,i_0}$. To prove the second one we use (iii)
to conclude that 
$ q_i(x) = q_i(x) p_{i_0}(x)  = (q_i p_{i_0})(x) = A(i,i_0) p_{i_0}(x) = c_{i_0}(i) =
\hat \rho_i (i_0, c_{i_0})$.

Suppose now that $p_i(x) = 0$ for every $i\in \Gen$, and define $b \in 2^{\Gen}$ by
$b(i) = q_i(x)$.
Clearly $p_i(x) = 0 = \hat\delta_i (\star,b)$ and
$q_i(x) = b(i) = \hat \rho_i(\star,b)$ for every $i\in \Gen$, so the only thing that needs
proving is that $(\star,b) $ is indeed in $ \tilde\Gamma_A$. For this it suffices to
 show that $b$ is an accumulation point of the collection of columns in $2^{\Gen}$.
Let 
$$
V = V(b,F) := \{ c\in 2^{\Gen}: c(j) = b(j) \text{ for } j\in F\}
$$
be a neighborhood of $b$ and define  
$F_0 := \{i\in F: b(i) =0\}$ and $F_1 :=  \{ i\in F: b(i) =1\}$,
so that the characteristic function of $V(b,F)$ is given by
$$
1_V(c) = \prod_{i\in F_1} c(i) \prod_{i\in F_0} (1-c(i))
$$
and we have that 
\stepcounter{theorem}\begin{equation}\label{contra}
1 = \prod_{i \in F_1} b(i) \prod_{i \in F_0} (1-b(i)).
\end{equation}

If there were only finitely many $j\in \Gen$ for which $c_j \in V$, so that the function
$$
 j \mapsto 1_V(c_j) = \prod_{i\in F_1} c_j(i) \prod_{i\in F_0} (1-c_j(i)) =
    \prod_{i\in F_1} A(i,j) \prod_{i\in F_0} (1-A(i,j)) 
$$
had finite support. 
Then relation (iv) would yield
$$
 \prod_{i\in F_1} q_i \prod_{i\in F_0} (1-q_i)  = 
 \sum_j  A(F_1, F_0, j)  p_j,
$$
and by evaluation at $x$ we would get 
$\prod_{i\in F_1} b(i) \prod_{F_0} (1-b(i)) = 0$, in contradiction with \eqref{contra}. Thus,
there are infinitely many $j\in \Gen$ for which $c_j \in V$, so that $(\star,b) $ is indeed in
$\tilde\Gamma_A$, which finishes the proof of the unital case.

In order to conclude that $\rowa$ is the universal, not necessarily unital, $C^*$-algebra
with the given presentation, it suffices to observe that if the point $(\star,\zerocol)$
is in $\tilde \Gamma_A$, then its image in $X$ is the unique point
$x_\infty$ defined by $p_i(x_\infty)= q_i(x_\infty)= 0$.
\end{proof}

\begin{proposition} \label{rels-revisit}
Let $A$ be a $0$-$1$ matrix with no identically zero rows.
The $C^*$-algebra $\oa$ is (canonically isomorphic to) the universal $C^*$-algebra generated by elements
$\{S_i: i\in \Gen\}$ such that the maps $\delta_i \mapsto S_i S_i^*$ and
$\rho_i \mapsto S_i^* S_i$ may be extended to a homomorphism of $\rowa$ into $C^*(\{S_i: i \in
\Gen\})$.
\end{proposition} 
\begin{proof}
Let $\{S_i: i\in \Gen\}$ be the generators of $\oa$. 
Then the projections $P_i = S_i S_i^*$ and $Q_i = S_i^* S_i$ satisfy the relations (i) --- (iv) of
\lemref{rowa-univ}, and hence there is a homomorphism of $\rowa$ into $\oa$ such that
  \stepcounter{theorem}\begin{equation}
  \label{canonical-map}
\delta_i \mapsto S_i S_i^* \quad \text{ and }\quad \rho_i \mapsto S_i^* S_i.
\end{equation}

Conversely, assume the $S_i$ are elements such that the maps in \eqref{canonical-map} 
extend to a homomorphism of $\rowa$.  
Then (ii) holds because the $\delta_i$ are pairwise  orthogonal projections,
which implies that the $S_i$ are partial isometries.  
Relation (i) holds because the $\rho_i$ commute with each
other, and  (iii) holds because $\rho_i \delta_j = A(i,j) \delta_j$. If the
finite sets $X$ and $Y$ of $\Gen$ are such that the function $A(X,Y,j)$ is finitely supported, 
then $\prod_X \rho_x \prod_Y (1 -\rho_y) = \sum_j A(X,Y,j) \delta_j$, and the 
existence of a homomorphism extending \eqref{canonical-map} implies that 
$$
\prod_{x\in X} S_x^*S_x \prod_{y\in Y} (\bfone - S_y^* S_y)\  =\ 
     \sum_{j\in\Gen} A(X,Y,j) S_j S_j^*,
$$
so (iv) also holds. By Theorem 8.6 of \cite{infinoa} the $S_i$ generate a 
homomorphic image of $\oa$.
\end{proof}

\section{The map $(I - A^t)$ and the $K$-theory of $\oa$.}\label{themap}

  Since the image of a basis vector under the linear transformation
  defined by a matrix 
  is the corresponding column of the matrix, and since columns of $A^t$
  correspond to rows of $A$, it makes sense to denote by $A^t$ 
  the  map  $\delta_i \mapsto \rho_i$.  We will define a map
  $(I - A^t)$ with domain  $\bigoplus_{i\in \Gen} {\mathbb Z}$ by simply saying that
  $(I - A^t) \delta_i = \delta_i - \rho_i $,
  but, of course, it is crucial to determine 
  the appropriate co-domain for such a map to be useful in computing the
  $K$-theory of $\oa$.

  \begin{definition}
   Let 
   $$
   \zrowa := {\operatorname{Ring}}\{\delta_i, \rho_i: i\in \Gen\} \subset \mathbb Z^\Gen
   $$
  be the ring of functions on $\Gen$ generated by the rows of the matrices $I$ and $A$, and let
  $$
   \zurowa := {\operatorname{Ring}}\{\onecol,\delta_i, \rho_i: i\in \Gen\} 
   $$
   be the corresponding unital ring.
   \end{definition}

  \begin{proposition} Under the isomorphism of $\urowa$ to $C(\tilde\Gamma_A)$
  given in \proref{SpectrumRowa}, we have that
  $  \zrowa \cong C_0(\Gamma_A;\mathbb Z) $ and 
   $ \zurowa \cong  C(\tilde\Gamma_A;\mathbb Z) $.
 \end{proposition}
 
 \begin{proof}
 Since  $\rowa = C_0(\Gamma_A)$ is the $C^*$-algebra generated by the projections
   $\delta_i$, and $\rho_i$, the result follows from the following elementary 
   general lemma.
   \end{proof}
   
  \begin{lemma}\label{subring}
  Let $\{E_i: i\in \Lambda\}$ be a family of commuting projections,
  let $A := C^*(\{E_i: i\in \Lambda\}$, and denote  by $X$ the spectrum of $A$.
  Then (identifying $C_0(X)$ to $A$ under the Gelfand  transform)  
  $C_0(X;\mathbb Z) = {\operatorname{Ring}}\{E_i: i\in \Lambda\}$.
  If the generating family is closed under multiplication,
  then $C_0(X;\mathbb Z) = \spa_{\mathbb Z}\{E_i: i\in \Lambda\}$.
  \end{lemma}
  
 \begin{proof}
 Let $\mathfrak E$ be the subring of $C_0(X)$ generated by the $E_i$; clearly $\mathfrak E
 \subset C_0(X;\mathbb Z)$.
 To prove that $\mathfrak E $ is all of $C_0(X;\mathbb Z)$ it is enough to show 
 that every projection $P\in C_0(X)$ lies in $\mathfrak E$.
 Let $P$ be a projection in $C_0(X)$ and let $\epsilon >0$.
 Then there exists a finite sum  $ \sum_{i\in I} \lambda_i F_i$, with $\lambda_i \in \mathbb C$ and
  each $F_i$ a product of $E_j$'s, such that $\|P - \sum_{i\in I}\lambda_i F_i\| < \epsilon $.
 For each $J\subset I$ let 
 $$
 F_J = \prod_{i\in J} F_i \prod_{i\in I\setminus J} (1 - F_i) \in M(A);
 $$
 then the $F_J$ are mutually orthogonal projections, and,  
 denoting the characteristic function of $J$ by $1_J$, we have 
 $F_i = \bigoplus_{J\in 2^I}1_J(i)  F_J $ for every $i\in \Lambda$.
 Thus 
 $$
\sum_{i\in I} \lambda_i F_i = \sum_{i,J} \lambda_i 1_J(i) F_J =
   \sum_{J\neq \emptyset}(\sum_i \lambda_i 1_J(i)) F_J = \sum_{J\neq \emptyset} \mu_J F_J,
   $$
   with $\mu_J := \sum_i \lambda_i 1_J(i)$.
   Since the $F_J$ are mutually orthogonal and $\sum_{J\neq \emptyset} \mu_J F_J$ is
   within $\epsilon$ of the projection $P$, 
   we have that $\mu_J \in (-\epsilon, \epsilon) \cup  (1-\epsilon, 1+\epsilon)$.
   Hence 
   $$
   \| P - \sum_{J: \mu_J \in (1-\epsilon, 1+\epsilon)}  F_J \| \leq 
   \|P - \sum_{J\neq \emptyset} \mu_J F_J \| + \| \sum_{J\neq \emptyset} \mu_J F_J
   - \sum_{J: \mu_J \in (1-\epsilon, 1+\epsilon)}  F_J \| < \epsilon + \epsilon .
   $$
   Since $P$ and $\sum_{J: \mu_J \in (1-\epsilon, 1+\epsilon)} F_J$
   are projections, it suffices to take $\epsilon < 1/2$ to conclude that they are equal,
   and hence that $P \in \mathfrak E$.
   
   The second assertion of the lemma follows easily because the subgroup generated
   by a multiplicatively closed family of projections is a subring.
 \end{proof}
  
  Next we consider the projections in $C(\ometil)$.
    For $X$ a finite subset of $\Gen$ let $Q_X: = \prod_{x\in X} Q_x$, 
    the empty product being equal to the identity by convention.
    We will consider elements of the form 
   $
   S_\mu Q_X S_\mu^*,
   $
   with $\mu \in \mathbb F^+$ an admissible path. 
   When $\mu \neq e$ we assume that $X$ contains the last letter
   $\mu_{|\mu|}$ of $\mu$; this has no effect on $ S_\mu Q_X S_\mu^*$ because
   $S_{\mu_{|\mu|}}^* S_{\mu_{|\mu|}}S_\mu^* = S_\mu^*$. It is clear that every
   $S_\mu Q_X S_\mu^*$ is a projection.
   
  \begin{proposition}  \label{SQS*} 
  The family 
  $\{S_\mu Q_X S_\mu^*: X\text{ a finite subset of }\Gen, \mu \in \mathbb F^+\} $
  is closed under multiplication,
  $C(\ometil) = \cspa_{\mathbb C}\{S_\mu Q_X S_\mu^*\}$, and
  $C(\ometil;\mathbb Z) = \spa_{\mathbb Z}\{S_\mu Q_X S_\mu^*\}$ 
  \end{proposition}
  \begin{proof} 
    Using  the relations (i), (ii) and (iii) from Section~\ref{prelim} one shows that 
   the product $(S_\mu Q_X S_\mu^*)( S_\nu Q_Y S_\nu^*) $
   vanishes unless $\mu \leq \nu $ or $\nu \leq \mu$. 

  Suppose $\mu \leq \nu$ and write
   $\nu = \mu \nu'$. Then the product is equal to $S_\mu Q_X S_{\mu}^* S_{\mu} S_{\nu'}Q_Y
  S_\nu^* =  S_\mu Q_X  S_{\nu'} Q_Y S_\nu^*$,
   because  $Q_x$ absorbs $S_{\mu}^* S_{\mu}$, which is equal to $S_{\mu_{|\mu|}}^*
  S_{\mu_{|\mu|}}$ with  $\mu_{|\mu|}$ the last letter of $\mu$.
  By relation (iii),
   $Q_X S_{\nu'}$ is either zero or $S_{\nu'}$. Hence 
   $(S_\mu Q_X S_\mu^*)( S_\nu Q_Y S_\nu^*) $ is either zero or $ S_\nu Q_Y S_\nu^* $.
   When $\nu \leq \mu $ we take adjoints and conclude that the product is
   either zero or $ S_\mu Q_X S_\mu^* $.
   
   The projections $P_i$ and $Q_i$ are in $\{S_\mu
   Q_X S_\mu^*\} $, so (ii) above holds,
  and (iii) follows by \lemref{subring}.
  \end{proof}
 
  \begin{theorem} 
  \label{ktheory-thm}
  Let $\Gen$ be a set and suppose $A= \{A(i,j)\}_{i,j\in \Gen}$
  is a $0$-$1$ matrix with no identically zero rows.
  Define a map 
  $(I - A^t) : \bigoplus_{i \in \Gen}{\mathbb Z} \to  \zrowa$
   by $$ (I - A^t) \delta_i = \delta_i - \rho_i, \qquad  i \in \Gen.
   $$
   Then 
  \begin{enumerate}
  
  \medskip
 \item[(i)] \ $ K_0(\oa) = \zrowa /\ima (I-A^t)$, and
 
 \medskip
 \item[(ii)] \  $ K_1(\oa) = \ker (I-A^t)$. 
 
 \medskip
  \end{enumerate}

  \end{theorem}
  
  \begin{proof}
  We will first work within the unital category and, viewing 
  $(I -A^t)$ as a map into $\zurowa$, we will show that
  \begin{eqnarray}\stepcounter{theorem}
  K_0(\uoa) &=& 
  \zurowa/\ima (I-A^t), \label{ukzero}\\
  \stepcounter{theorem}
  K_1(\uoa) &=& \ker (I-A^t). \label{ukone}
  \end{eqnarray}
  This will prove (ii) because $K_1(\oa)  = K_1(\uoa)$. 
  To prove (i) recall that when $\oa \neq \uoa$, we have that
  $K_0(\oa) := \ker (\varepsilon_*: K_0(\uoa) \to K_0(\mathbb C))$. We will see that the
  homomorphism of $\zurowa/\ima (I-A^t)$ to $\mathbb Z$ corresponding to
  $\varepsilon_*$ sends $\delta_i$ and $\rho_i$ to zero for every $i\in \Gen$, hence
  its kernel is $\zrowa/\ima (I-A^t)$.

  As before, we will denote by $L_\alpha$ the map
  $$
  \sum_{x\in {\Gen}}(\idx - \alpha_x\inv): 
  \bigoplus_{x\in \Gen} C(\Delta_x,{\mathbb Z})  \to  C(\ometil,{\mathbb Z}).
  $$
  {}From  \proref{pvseq-lim} we know that 
  $K_0(\uoa) = C(\ometil,{\mathbb Z}) / \ima L_\alpha$ and
  $K_1(\uoa) = \ker L_\alpha$. In order to compute these in terms of $A$ we use
  the following diagram (cf. \cite[p.33]{cun2}):
  \stepcounter{theorem}\begin{equation}
  \label{sq-diagram}
  \begin{CD}
  \bigoplus_{x\in \Gen} C(\Delta_x ;{\mathbb Z}) 
  @>{\ \ L_\alpha \ \ }>> C(\ometil;{\mathbb Z}) \\
  @A{\jay}AA  @A{\psi}AA\\
  \bigoplus_{x\in \Gen} {\mathbb Z} @>{ I-A^t }>>  \zurowa.  
  \end{CD}
  \end{equation}
  
  \medskip
  The vertical arrows are from \proref{psihat}, and are determined by $\jay (\delta_x) = P_x$,
  $\psi (\delta_x) = P_x$
   and
  $\psi (\rho_x) = Q_x$. The diagram commutes because 
  $\psi (I-A^t) \delta_x = \psi(\delta_x - \rho_x) = P_x - Q_x$ and
  $ L_\alpha \jay \delta_x = L_\alpha P_x = P_x - \alpha_x\inv(P_x) = P_x - Q_x$.
  That the kernel and co-kernel of $ L_\alpha$ are isomorphic to
  those of $(I-A^t)$ will follow from the following general 
  lemma about commuting diagrams of group homomorphisms (we omit the straightforward proof).
  
 \begin{lemma} Let  
  \begin{equation*}
  \begin{CD}
 G_1
  @>{\ \ L \ \ }>> G_2\\
  @A{\jay}AA  @A{\psi}AA\\
  H_1 @>{\ \  M \ \ }>>  H_2
  \end{CD}
  \end{equation*}
  be a commuting diagram of abelian groups in which the vertical arrows are injective and such that
  \begin{enumerate}
  \item[(I)] every $x\in G_2$ is equivalent, modulo $L(G_1)$, to an element of $\psi(H_2)$, and
  \item[(II)]  for all $x\in G_1$, if $L(x) \in \ima (\psi)$ then $x\in \ima(\jay)$.
  \end{enumerate}
  Then $\ker M = \ker L$ and $\coker M = \coker L$.
  \end{lemma}
  
  The proof of the theorem is finished by verifying, in the next two lemmas, that the diagram
\eqref{sq-diagram}
  satisfies the hypothesis (I) and (II) above.
  \end{proof}
  
  \begin{lemma}\label{moduloL}
  Every element of $C(\ometil, {\mathbb Z})$ is equivalent, 
  modulo $ L_\alpha \left(\bigoplus_x C(\Delta_x,{\mathbb Z})\right)$,
  to an element of $\psi(\zurowa)$.
  \end{lemma}

\begin{proof}
  By \proref{SQS*} it suffices to prove the statement for elements of the form  $S_\mu Q_X S_\mu^*$.

  Suppose $\mu \neq e$ and write $\mu = \mu_1 \nu$
  with $\mu_1$ in $\Gen$.
  It is easy to see that $S_\mu Q_X S_\mu^*$
  is in the domain of $\alpha_{\mu_1}\inv$ and that
 $\alpha_{\mu_1}\inv (S_\mu Q_X S_\mu^*) = S_\nu Q_X S_\nu^*$,
   with $\nu = \mu_2 \mu_3 \cdots \mu_{|\mu|}$, so that 
  $$
S_\mu Q_X S_\mu^* - S_\nu Q_X S_\nu^* =  L_\alpha ( S_\mu Q_X S_\mu^* ).
$$ 
  We may now continue the process, removing one letter of $\mu$ at each step.
  When removing the last letter of $\mu$ we need to use our assumption that
  it is in $X$; this ensures that $Q_X$ is in the image 
   of $\alpha_{\mu_{|\mu|}\inv}$ so that the last step gives
   $\alpha_{\mu_{|\mu|}}\inv (S_{\mu_{|\mu|}} Q_X S_{\mu_{|\mu|}}^*) =  Q_X 
   $, and hence
  $$
  S_{\mu_{|\mu|}} Q_X S_{\mu_{|\mu|}}^* -  Q_X  =  L_\alpha ( S_{\mu_{|\mu|}} Q_X S_{\mu_{|\mu|}}^* ).
  $$ 
  We can then use an elementary ``telescopic sum" argument to conclude that
  $$
  S_\mu Q_X S_\mu^* -  Q_X \in 
  L_\alpha \left( \textstyle{\bigoplus_x} C(\Delta_x;\mathbb Z) \right),
  $$
  which finishes the proof because 
  $Q_X  = \prod_{x\in X} Q_x = \psi (\prod_{x\in X} \rho_x) \in \psi(\zurowa)$.
  \end{proof}  

  \begin{lemma}\label{LfinRA}
  Let $ f \in \bigoplus_x C(\Delta_x;\mathbb Z)$ be such that  $ L_\alpha f \in
  \psi(\zurowa)$. Then $f \in \jay(\bigoplus_{x\in \Gen}\mathbb Z)$.
  \end{lemma}

  \begin{proof}
  We view $f = (f_x)_{x\in \Gen}$ as a $\mathbb Z$-valued function on $\ometil$ in the obvious way.
  Assume that $f\notin \spa\{ P_x: x\in {\Gen}\}$.  We claim that
  there exist $n\geq 1$
  and a pair of unbounded elements $\xi $ and $\eta$ of $\ometil$ such that:
  \begin{enumerate}
  \item[(i)]  $\sigma(\xi)|_n = \sigma(\eta)|_n$,
 \item[(ii)] $f(\xi) \neq f(\eta)$, and
 \item[(iii)] $n$ is maximal in the sense that if $\xi'$ and $\eta'$ are two unbounded
  elements of $\ometil$ such that $\sigma(\xi')|_{n+1} = \sigma(\eta')|_{n+1}$, then
  $f(\xi') = f(\eta')$.
 \end{enumerate}
  To prove the claim we argue as follows.
  By assumption $f$ is not constant in at least one of the
  $\Delta_x$; since unbounded elements are dense, 
  there exist unbounded elements $\xi$ and $\eta$ in the same $\Delta_x$, 
  such that  (ii) holds and
  (i) holds for $n = 1$ because $\sigma(\xi)|_1 = x = \sigma(\eta)|_1$. 
  To see that there is a maximal such $n$ it suffices to show
   that the collection of all the $n\in \mathbb Z$ for which there exist
  $\xi$ and $\eta$  satisfying (i) and (ii) is bounded. 
   By \proref{SQS*}, the function $f$  is a finite sum of the
  form  $f = \sum S_\mu Q_X S_\mu^*$. 
  If $\sigma(\xi')|_{m} = \sigma(\eta')|_{m}$ for $m$ strictly larger than 
  the lengths of all the $\mu$'s
  appearing in the sum, 
 then $(S_\mu Q_X S_\mu^*)(\xi') = (S_\mu Q_X S_\mu^*)(\eta')$ from which 
  it follows that $f(\xi') = f(\eta')$. This concludes the proof of the claim.

  By definition for $f = (f_x)_{x\in \Gen}$ we have
  $$
   L_\alpha f = \sum_{x\in \Gen} (\idx -\alpha_x\inv) f_x =
   \sum_{x\in \Gen} f_x - \sum_{x\in \Gen} \alpha_x\inv(f_x),
  $$ 
  where the set $F: = \{x \in \Gen : f_x \neq 0\}$  is
  finite because $(f_x) \in \bigoplus_x C(\Delta_x,\mathbb Z)$.
  Evaluation at $\xi$ gives
  $$
  L_\alpha f (\xi) = 
               f(\xi) - \sum_{\substack{x\in F \\  x\inv \in \xi}}
                       f(x\xi), 
  $$
  so
 \stepcounter{theorem}\begin{equation}
 \label{rlxi}
  f(\xi)=\sum_{\substack{x\in F \\  x\inv \in \xi}}
                f(x\xi) + L_\alpha f (\xi),
  \end{equation}
  and similarly,
  \stepcounter{theorem}\begin{equation}
  \label{rleta}
  f(\eta)=
  \sum_{\substack{x\in F \\  x\inv \in \eta}} 
                f(x\eta) +L_\alpha f (\eta).
  \end{equation}
  The sums for $\xi$ and for $\eta$ are over the same finite sets because 
  $x\inv \in \xi $ if and only if $A(x,\sigma(\xi)|_1) =1$ and we know that 
  $\sigma(\xi)|_1 = \sigma(\eta)|_1$.
  
   By hypothesis, $ L_\alpha f$ is in $\psi (\zurowa) $, which is the unital subring of 
   $C(\ometil; \mathbb Z)$ generated by the projections $\{P_i, Q_i: i \in \Gen\}$.
   Since $P_i(\xi) = P_i(\eta) $ and $Q_i(\xi) = Q_i(\eta)$, we have that
   $L_\alpha f(\xi) = L_\alpha f(\eta)$. Moreover, the stems of $x\xi$ and $x\eta$
   coincide up to the place $n+1$ because
  $$
  \sigma(x\xi)|_{n+1} = x\sigma(\xi)|_n 
  =  x\sigma(\eta)|_n= \sigma(x\eta)|_{n+1},
  $$
   so property  (iii) implies that $f(x\xi)= f(x\eta)$. 
   But then the right hand sides of \eqref{rlxi} and \eqref{rleta} 
  coincide, and we must have
  $f(\xi) = f(\eta)$, which contradicts property (ii). 
  Thus no such $\xi$ and $\eta$ exist, so 
  $f$ is constant on the subset $\Delta_x$ for every $x\in {\Gen}$.
  \end{proof}

  \section{Examples with edge matrices}

  In this section we consider a particular case of matrices for which the computations are
  relatively easy to carry out.  Specifically, we assume $A$ to be an {\em edge matrix}.
  By definition, the index set of an edge matrix $A$ 
  is the edge set of a directed graph, with $A(i,j) = 1$ if 
  the range of $i$ is equal to the source of $j$, and $A(i,j) =0$ otherwise.
  As a consequence, any two rows of $A$
  are either equal (when they correspond to edges with the same range),
  or else orthogonal (when they correspond to edges with different ranges).   
  Because of this feature, the $C^*$-algebra $\rowa$, the ring $\zrowa$, and hence the 
  $K$-theory of $\oa$ have a rather simple form. We begin with a general observation,
  and then restrict our attention to some examples.

 \begin{proposition}
 \label{RowaEdge}
 Suppose $A$ is an edge matrix over $\Gen$.
 Then $\rowa = \cspa \{\rho_i, \delta_i: i \in {\Gen}\}$,
  and $\zrowa = \spa_{\mathbb Z}\{\rho_i, \delta_i: i \in {\Gen}\}$
  (the group generated by  $\rho_i$'s and $\delta_i$'s).
  \end{proposition}
\begin{proof}
Since the rows of $A$ are either equal or orthogonal, 
the collection of $\rho_i$'s and $\delta_i$'s is closed
under products.
\end{proof}
  
  The only non identically zero accumulation points
  of the collection of columns are the columns that appear
  infinitely many times; it is easy to see that they 
  correspond to the infinite rows. 
  In this case, $\Gamma_A$ is
  obtained by adding a point $(\star,c)$ to $\{(i,c_i): i\in \Gen\}$ for each column $c$ that
  appears infinitely many times in $A$. 
  We compute next $K_0(\oa)$ and $K_1(\oa)$ for a few edge matrices.
  
  \begin{example}
   Let $A$ be the $\N \times \N$ matrix all of whose entries are $1$. By \proref{RowaEdge}
   we have that $\zrowa = \bigoplus_{\N} \mathbb Z + \mathbb Z\cdot \onecol$.
   Let $f\in \bigoplus_{\N} \mathbb Z$.
   Then $(I-A^t) f = f - (\sum_{\mathbb N} f(n) ) \onecol$.
   It follows that $f$ is in $ \ker (I-A^t)$ if
   and only if it is constant, and hence identically zero. 
   Suppose now $h \in \zrowa$ and denote the coefficient of $\onecol$
   by $h_\infty$, so that $f: = h - h_\infty \onecol$ is finitely supported. 
   If $h_\infty = -\sum_{\N} (h - h_\infty \onecol) (n)$ then
   $h = f - (\sum_{\N} f(n) ) \onecol = (I - A^t) f$. 
   Conversely, if  $h = f - (\sum_{\N} f(n) ) \onecol$
   for some $f\in \bigoplus_{\N} \mathbb Z$, then 
   $h_\infty = - \sum_{\N} f(n)$ and $f = h - h_\infty \onecol$, so that 
 $h_\infty = -\sum_{\N} (h - h_\infty \onecol)$.  
   Thus  $\ima (I-A^t)$ is the kernel of the homomorphism
   $\varphi:  \zrowa \to  \mathbb Z$ given by
   $\varphi(h) =   h_\infty +\sum_{\N} (h - h_\infty \onecol)(n) $, from which it follows
  that $\zrowa / \ima (I-A^t) = \mathbb Z$.
  
  Since in this case  $\mathcal O_{A} = {\mathcal O}_\infty$, the above computation
  recovers the classical result, from \cite{cun2}, that $K_0(\mathcal O_\infty) = \mathbb Z$
   and $K_1(\mathcal O_\infty) = 0$.
  \end{example}
   
  \begin{example} Consider now the following irreducible ``checkerboard" matrix 
  of alternating zeros and ones.
  $$ 
  A = 
  \left( 
  \begin{array}{ccccccc}
  0&1&0&1&.&.&.\\
  1&0&1&0&.&.&.\\
  0&1&0&1&.&.&.\\
  1&0&1&0&.&.&.\\
  .&.&.&.&.& & \\
  .&.&.&.& &.& \\
  .&.&.&.& & &.
  \end{array}
  \right)
  $$
  We will show that 
  $K_0(\oa) = \mathbb Z \times \mathbb Z$ and $K_1(\oa)= 0$.
  
  Denote by $\onecol_{\even}$ (respectively, $\onecol_{\odd}$)
  the characteristic function of the 
  set of even numbers (respectively, the set of odd numbers).
  By \proref{RowaEdge} we have that
  $$
  \zrowa = \textstyle\bigoplus_{\N} {\mathbb Z} + \mathbb Z \cdot{\onecol}_{\odd} +
  \mathbb Z \cdot{\onecol}_{\even}.
  $$
  By definition, for $f\in \bigoplus_{\N} {\mathbb Z}$,
  $$
  (I-A^t) f =
  f - \textstyle\sum_{n \in \odd} f(n) \cdot \onecol_{ \even} - 
  \textstyle\sum_{n \in  \even} f(n) \cdot \onecol_{\odd}.
  $$

  Since $f$ is finitely supported, if  $(I-A^t)f = \zerocol $ 
  then $ \textstyle\sum_{\odd} f(n) =  \textstyle\sum_{\even} f(n) = 0$, and hence $f$ itself
is zero, from which it follows that  $K_1(\oa) = 0$.
  
 Let  $ h \in \zrowa$, and, for convenience of notation, write 
  $h = \sum_{n\in \N} \lambda_n  \delta_n - \mu_e \onecol_{\even} - \mu_o
    \onecol_{\odd}$, so that  $h(n) = \lambda_n  - \mu_e$ for $n$ even,
    $h(n) = \lambda_n - \mu_o$ for $n$ odd, and the $\lambda$'s are uniquely
    determined by the requirement that the
    sequence $\lambda_n$ be finitely supported. 
    If $h$ is in the image of $(I-A^t)$, then there exists $f \in \bigoplus_{\N}
    {\mathbb Z}$ such that 
   $$ 
   f - \textstyle\sum_{\odd} f(n) \cdot \onecol_{\even} - \textstyle\sum_{\even} f(n) \cdot
   \onecol_{\odd} 
      \    = \ 
   \sum_{n\in \N} \lambda_n \delta_n - \mu_e
    \onecol_{\even} - \mu_o \onecol_{\odd}.
    $$
    It follows that $f = \sum_{n\in \N} \lambda_n \delta_n $ and that 
      $\textstyle\sum_{\odd} f(n) = \mu_e$ and $\textstyle\sum_{\even} f(n) =\mu_o$.
    It is now easy to see that the image of $(I-A^t)$ is the kernel 
      of the homomorphism 
      $\varphi : \zrowa \to \mathbb Z \times \mathbb Z$ given by
      $$
     \varphi_1(h) =  \mu_e - \textstyle\sum_{\odd} \lambda_n  
     $$
     and by
     $$
     \varphi_2(h) =   \mu_o - \textstyle\sum_{\even} \lambda_n.
     $$
    Since $\varphi$ is obviously surjective,  $K_0(\oa) = \mathbb Z \times \mathbb Z$.
    \end{example}
  
  \begin{example}
  Interchanging zeros and ones in the preceding matrix 
  gives another checkerboard matrix,
  $$ 
  A = 
  \left( 
  \begin{array}{ccccccc}
  1&0&1&0&.&.&.\\
  0&1&0&1&.&.&.\\
  1&0&1&0&.&.&.\\
  0&1&0&1&.&.&.\\
  .&.&.&.&.& & \\
  .&.&.&.& &.& \\
  .&.&.&.& & &.
  \end{array}
  \right).
  $$ 
  One could carry out the same sort of analysis as in the previous case, but it is easier to 
  observe that the index set may be reshuffled by listing first all
  the odd numbers and then all the even ones,
  in such a way that the rearranged matrix is the  
   direct sum of two matrices of `all ones'. From this we conclude that
  $\oa$ is  the direct sum of two copies of 
  $\mathcal O_\infty$, so that $K_0(\oa) = \mathbb Z \times \mathbb Z$ and $K_1(\oa) = 0$.

Notice that this $\oa$ is not simple, while the preceding one is simple
and  purely infinite, by \cite[Theorems 14.1 and 16.2]{infinoa}.
  \end{example}
  
 \begin{remark} In view of the simplifications in the calculations involving
  edge matrices, it is natural to ask
  whether every $\oa$ is the Cuntz-Krieger algebra of a
  conveniently chosen edge matrix. 
  It is shown in Proposition 4.1 of \cite{mrs} (see also \cite{ror}),
  that given an arbitrary {\em finite} matrix $A$ of zeros and ones
  one may find an associated 
  edge matrix $B$ (also finite, but in general of a different size)
  with ${\mathcal O}_B \cong \oa$. However, the method used to prove this
  equivalence breaks down for  
  infinite matrices that are not row-finite, because the formulas
  for such matrices involve nonconvergent infinite sums.

\section*{Acknowledgment}
This work was completed during a visit of M.L. to Florian\'{o}polis, and we thank 
the Department of Mathematics at Newcastle for partially supporting this visit.
  
  \end{remark}

\end{document}